\documentclass[12pt,reqno]{amsart}
\usepackage{amssymb}
\usepackage{amsfonts}
\usepackage{amsmath}
\usepackage{amsthm}
\usepackage{bm}
\usepackage{color}
\usepackage{enumerate}
\usepackage{epsfig}
\usepackage{hyperref}
\usepackage{bm}
\makeatletter

\newcommand{\Aut}{\mathrm{Aut}}
\newcommand{\Cay}{\mathrm{Cay}}

\newcommand{\BCay}{\mathrm{BCay}}
\newcommand{\GL}{\mathrm{GL}}

\newtheorem{thm}{Theorem}[section]
\newtheorem{cor}[thm]{Corollary}
\newtheorem{pro}[thm]{Proposition}
\newtheorem{lem}[thm]{Lemma}

\newtheorem{rem}[thm]{Remark}
\begin{document}

\title[Isomorphisms of bi-Cayley graphs on generalized quaternion groups]
{Isomorphisms of bi-Cayley graphs on generalized quaternion groups}

\author{Jin-Hua Xie}
\address{Jin-Hua Xie, Center for Combinatorics, LPMC, Nankai University, Tianjin 300071, China}
\email{jinhuaxie@nankai.edu.cn (J.-H. Xie)}


\author{Zhishuo Zhang}
\address{Zhishuo Zhang, School of Mathematics and Statistics, The University of Melbourne, Parkville, VIC 3010, Australia}
\email{zhishuoz@student.unimelb.edu.au (Z. Zhang)}


\begin{abstract}
Let $G$ be a finite group and $S$ be a subset of $G$. The bi-Cayley graph $\BCay(G,S)$ is the graph with vertex set $G\times \{0,1\}$ and edge set $\{\{(x,0),(sx,1)\}\mid x\in G,s\in S\}$. A bi-Cayley graph $\BCay(G,S)$ is called a BCI-graph if for every $T\subseteq G$, the isomorphism $\BCay(G,S)\cong \BCay(G,T)$ implies that $T=gS^\alpha$ for some $g\in G$ and $\alpha\in \Aut(G)$. We say a group $G$ an $m$-BCI-group if every bi-Cayley graphs of $G$ with valency at most $m$ is a BCI-graph. In this paper, we show that for $m\in\{2,3\}$, the generalized quaternion group of order $4n$ with $n\geq 2$ is an $m$-BCI-group if and only if it is an $m$-DCI-group if and only if it is an $m$-CI-group if and only if $n$ is odd or $n=2$.


\vskip 5pt

\end{abstract}
\maketitle

\section{Introduction}

Throughout this paper, all digraphs and graphs are finite, with no parallel edges or loops, and all groups are finite. For a positive integer $m$, a digraph is said to have \emph{valency $m$} if every vertex has both in-valency and out-valency equal to $m$. Let $G$ be a group and let $S$ be a subset of $G$ with $1\notin S$. The \emph{Cayley digraph} of $G$ with respect to $S$, denoted by $\Cay(G,S)$, is a digraph with vertex set $G$ and arc set $\{(g,sg)\mid g\in G,~s\in S\}$. If $S=S^{-1}$, then $(u,v)$ is an arc if and only if $(v,u)$ is an arc. In this case, the Cayley digraph $\Cay(G,S)$ is a graph, called \emph{Cayley graph}.

We say that a Cayley digraph $\Cay(G,S)$ is a \emph{CI-digraph} if for every $T\subseteq G$, the isomorphism $\Cay(G,S)\cong \Cay(G,T)$ implies that there exists $\sigma\in\Aut(G)$ such that $S^{\sigma} = T$. The group  $G$ is called a \emph{DCI-group} if every Cayley digraph of $G$ is a CI-digraph. For a positive integer $m$, we refer to $G$ as an \emph{$m$-DCI-group} if every Cayley digraph of valency at most $m$ is a CI-digraph. We define \emph{CI-graph}, \emph{CI-group} and \emph{$m$-CI-group} in exactly the same way when $\Cay(G,S)$ is a Cayley graph.

The problem of determining which finite groups are DCI-groups or CI-groups has been extensively studied over the past half century. In 1967, \'Ad\'am~\cite{Ad} conjectured that every finite cyclic group is a DCI-group. Although this conjecture was disproven by Elspas and Turner~\cite{Elspas} three years later, it sparked interest among researchers in identifying which finite groups are DCI-groups or CI-groups (see~\cite{Alspach,Babai,C.H.Li1} for example). This problem was ultimately solved by Muzychuk~\cite{Mu1,Mu2} in 1997. He provided a complete classification of cyclic DCI-groups and CI-groups, stating that a cyclic group $C_n$ of order $n$ is a CI-group if and only if it is a DCI-group or $n\in\{8,9,18\}$, and is a DCI-group if and only if $n/\gcd(2,n)$ is square-free. Inspired by this, classifying finite DCI-groups and CI-groups has become a hot topic, see~\cite{DM,DE1,C.H.Li8,LLP} for example. However, it remains highly challenging, even for elementary abelian groups (see~\cite{Feng,Somlai,Spiga}). For more research on DCI-groups or CI-groups, we refer the reader to~\cite{DE0,KMSRS,XFRL,XFK,XFZ,XFK1}.

Similar to the construction of Cayley graphs, another interesting family of graphs is bi-Cayley graphs, on which isomorphism problems are also actively studied over the past 15 years. Let $G$ be a group and let $S$ be a subset of $G$ (possibly including the identity element). The \emph{bi-Cayley graph} of $G$ with respect to $S$, denoted by $\BCay(G,S)$, is a graph with vertex set $G\times \{0,1\}$ and edge set $\{\{(x,0),(sx,1)\}\mid x\in G,s\in S\}$. A bi-Cayley graph $\BCay(G,S)$ is called a \emph{BCI-graph} if for every $T\subseteq G$, the isomorphism $\BCay(G,S)\cong \BCay(G,T)$ implies that $T=gS^\alpha$ for some $g\in G$ and $\alpha\in \Aut(G)$. We say a group $G$ a \emph{BCI-group} if every bi-Cayley graphs of $G$ is a BCI-graph, and we call $G$ an \emph{$m$-BCI-group} if every bi-Cayley graphs of $G$ with valency at most $m$ is a BCI-graph.

The study of BCI-groups was initiated in 2008 by the authors of~\cite{XJSZL}, who proved that all finite groups are $1$-BCI-groups and provided a necessary and sufficient condition for a finite group to be a $2$-BCI-group. Later in 2010, Jin and Liu~\cite{JinL} demostrated that $\mathrm{A}_5$ is the only $3$-BCI-group among nonabelian simple groups. In 2019, Koike and Kov\'{a}cs~\cite{Koike} showed that a finite nilpotent group is a $3$-BCI-group if and only if it is in the form $U \times V$, where $U$ is a direct product of cyclic groups of the same odd order and $V\in \{1, C_{2^r}, C_2^r,\mathrm{Q}_8\}$. For more results concerning $m$-BCI-groups, we refer the reader to~\cite{JinL1,KK,LuWX}.

We are interested in the relationship between the isomorphism problems on Cayley digraphs (graphs) and bi-Cayley graphs of the same group. It is evident that a CI-group may not necessarily be a BCI-group. In fact, the cyclic group $C_8$ is a CI-group~\cite[Theorem 8.9]{C.H.Li8} but not a $4$-BCI group~\cite[Remark 1.2]{KK}. To our knowledge, every known $m$-BCI-group is both a DCI-group and a CI-group. In 2015, Arezoomand and Taeri~\cite{AT} formally conjectured that every BCI-group is a CI-group. This conjecture is still open\footnote{In the unpublished paper~\cite{SAMI}, Theorem~1 claims that the conjecture is true. However, there is an obvious mistake in the last two lines of the proof, which we are unable to correct.} and not many relevant results have been publish until now. In~\cite{ABGR,Qu}, it is shown that for each $m\in \{2,3\}$, the dihedral group of order $2n$ with $n\geq 3$ is an $m$-BCI-group if and only if it is a $m$-DCI-group if and only if it is a $m$-CI-group if and only $n$ is odd. In this paper, we prove a similar result for generalized quaternion groups (defined below).

For an integer $n$ with $n\geq 2$, define the generalized quaternion group (also called the dicyclic group) of order $4n$ by
\begin{equation}\label{eq:Q4n}
\mathrm{Q}_{4n}=\langle a,b\mid a^{2n}=1,\,b^2=a^n,\,a^b=a^{-1}\rangle.
\end{equation}
In particular, $\mathrm{Q}_8$ is the quaternion group. Note by~\cite[Theorem~1.1]{SG} that $\mathrm{Q}_8$ is a DCI-group and by~\cite[Lemmas~3.2 and~3.4, Corollary~3.5]{Ma} that for $n\geq 3$ and $m\in\{2,3\}$, $\mathrm{Q}_{4n}$ is an $m$-DCI-group if and only if it is an $m$-CI-group if and only if $n$ is odd\footnote{The proof of necessity of~\cite[Corollary~3.5]{Ma} is not shown in that paper. Although it is not difficult, we add it here for completeness. It suffices to show that when $n$ is even, $\mathrm{Q_{4n}}$ is not a $2$-CI-group. Clearly, $\Cay(\mathrm{Q_{4n}},\{a^{n/2},a^{-n/2}\})\cong \Cay(\mathrm{Q}_{4n},\{b,b^{-1}\})$. However, $\{a^{n/2},a^{-n/2}\}^\alpha\neq \{b,b^{-1}\}$ for any $\alpha\in\Aut(\mathrm{Q}_{4n})$ (see Lemma~\ref{lm:Q4n}\eqref{enu:fQ4n}). Therefore, $\mathrm{Q_{4n}}$ is not a $2$-CI-group.}. The following theorem determines which generalized quaternion group is an $m$-BCI-group, where $m\in\{2,3\}$.

\begin{thm}\label{mainth1}
Let $\mathrm{Q}_{4n}$ be as in~\eqref{eq:Q4n} for some integer $n\geq 2$. Then $\mathrm{Q}_{4n}$ is a 3-BCI-group if and only if it is a $2$-BCI-group if and only if $n$ is odd or $n=2$.
\end{thm}

Therefore, with all the above results, we have the following corollary.

\begin{cor}\label{cor1}
Let $\mathrm{Q}_{4n}$ be as in~\eqref{eq:Q4n} for some integer $n\geq 2$. Then for $m\in\{2,3\}$, the group $\mathrm{Q}_{4n}$ is an $m$-BCI-group if and only if it is an $m$-DCI-group if and only if it is an $m$-CI-group if and only if $n$ is odd or $n=2$.
\end{cor}

In Section~\ref{sec:2}, we introduce some preliminary results needed in this paper. In Section~\ref{sec:3}, we demonstrate that two particular bi-Cayley graphs are not isomorphic, in preparation for the proof of Theorem~\ref{mainth1} in Section~\ref{sec:4}.

\section{Preliminaries}\label{sec:2}
In this section we collect some basic facts that will be used later.
Throughout this paper, denote by $\mathbb{Z}_n$ the additive group of integers modulo $n$ and $\mathbb{Z}_n^*$ the multiplicative group of all integers coprime to $n$ in $\mathbb{Z}_n$. For a group $G$ and an element $a\in G$, we use $|a|$ and $|G|$ to represent the order of $a$ and $G$, respectively.

We first recall an elementary result on cyclic groups.
\begin{lem}\label{lm:cyclic}
  Let $G=\langle c\rangle$ be a cyclic group of order $m$ for some positive integer $m$. Then the following statements hold.
  \begin{enumerate}[\rm(i)]
    \item Every subgroup of $G$ is cyclic.
    \item $G$ has a unique subgroup of order $d$ for each divisor $d$ of $m$.
  \end{enumerate}
\end{lem}

The next lemma collects some properties on the generalized quaternion group $\mathrm{Q}_{4n}$. For brevity, we will reference the lemma only when a reminder may benefit the reader.

\begin{lem}\label{lm:Q4n}
  Let $\mathrm{Q}_{4n}$ be as in~\eqref{eq:Q4n} for some integer $n\geq 2$. Then the following hold.
  \begin{enumerate}[\rm(i)]
    \item \label{enu:aQ4n} $\mathrm{Q}_{4n}=\{b^ka^j\mid k\in\mathbb{Z}_2,~j\in \mathbb{Z}_{2n}\}$.
    \item \label{enu:bQ4n} $|a^j|=2n/\gcd(j,2n)$ and $|ba^j|=4$ for each  $j\in\mathbb{Z}_{2n}$.
    \item \label{enu:cQ4n} If $n\geq 3$, then $\langle a\rangle$ is characteristic in $\mathrm{Q}_{4n}$.
    \item \label{enu:dQ4n} If $n$ is odd and $k\in\mathbb{Z}_{2n}\setminus\{0,n\}$, then the subgroup $\langle a^k,b\rangle$ of $\mathrm{Q}_{4n}$ is isomorphic to $\mathrm{Q}_{4m}$, where
    \[
      m=\begin{cases}
        |a^k|/2 &\text{if }\, 2\nmid \gcd(k,2n)\\
        |a^k|&\text{if }\, 2\mid \gcd(k,2n).
      \end{cases}
    \]
    \item \label{enu:eQ4n} If $n$ is odd, then every isomorphism $\sigma$ between two subgroups of $\mathrm{Q}_{4n}$ can be extended to an automorphism of $\mathrm{Q}_{4n}$.
    \item \label{enu:fQ4n} If $n\geq 3$, then $\mathrm{Aut}(\mathrm{Q}_{4n})=\{\sigma_{r,s}\mid r\in \mathbb{Z}_{2n}^*,~s\in \mathbb{Z}_{{2n}}\}$, where $\sigma_{r,s}: a^j\mapsto a^{rj},~ba^j\mapsto ba^{rj+s}$.
  \end{enumerate}
\end{lem}
\begin{proof}
  Parts~\eqref{enu:aQ4n}--\eqref{enu:cQ4n} are well-known. Part~\eqref{enu:dQ4n} can be easily verified by~\eqref{eq:Q4n}, Part~\eqref{enu:eQ4n} is shown in~\cite[Lemma~2.4]{Ma}, and Part~\eqref{enu:fQ4n} can be seen from~\cite[Lemma~3.1]{XFK}.
\end{proof}


Two elements $x$ and $y$ of $G$ are said to be \emph{fused} if $x=y^\sigma$ for some $\sigma\in \Aut(G)$ and to be \emph{inverse-fused} if $x=(y^{-1})^\sigma$ for some $\sigma\in \Aut(G)$. A group $G$ is called \emph{FIF-group} if for arbitrary two elements of the same order are fused or inverse-fused. The following proposition gathers some results from \cite[Proposition 2.7]{ABGR} and \cite[Remark 1.2]{KK}.

\begin{pro}\label{finitegroup-1-2-3-BCI}
  The following statements hold.
  \begin{enumerate}[{\rm (i)}]
  \item \label{enu:1prop} Every finite group is $1$-BCI-group.
  \item \label{enu:2prop} A finite group is a $2$-BCI-group if and only if it is a FIF-group.
  \item \label{enu:3prop} Every finite cyclic group $C_n$ is a $3$-BCI-group.
  \item \label{enu:4prop} A finite cyclic group $C_n$ is a $4$-BCI-group if and only if $n$ is indivisible by $8$.
  \end{enumerate}
\end{pro}

The next lemma follows from~\cite[Lemmas~2.1(3) and~2.2]{LuWX} and~\cite[Lemma~2.8]{JinL}.

\begin{lem}\label{isomorphism-Bcay}
Let $G$ be a group. Then the following statements hold.
\begin{enumerate}[\rm(i)]
  \item \label{enu:aBiCay} For $S\subseteq G$, we have $\BCay(G, S)\cong \BCay(G,gS^\alpha)$ for all $\alpha\in \Aut(G)$ and all $g\in G$.
  \item \label{enu:bBiCay} If $S$ and $T$ are two subsets of $G$ containing the identity $1$, then $\BCay(G,S)\cong \BCay(G,T)$ if and only if $\BCay(\langle S\rangle, S)\cong \BCay(\langle T\rangle,T)$.
\end{enumerate}
\end{lem}

\begin{rem}\label{rmk:2.7}
  Lemma~\ref{isomorphism-Bcay}\eqref{enu:aBiCay} enables us to assume that $1\in S$ in the study of isomorphisms of $\BCay(G,S)$.
\end{rem}

\section{Non-isomorphic bi-Cayley graphs}\label{sec:3}

In this section, we aim to prove the following proposition by considering the spectrums (defined below) of graphs.
\begin{pro}\label{prop:not iso}
  Let $\mathrm{Q}_{4n}$ be as in~\eqref{eq:Q4n} for some integer $n\geq 2$. Then the two bi-Cayley graphs $\BCay(\mathrm{Q}_{4n},\{1,a,b\})$ and $\BCay(\mathrm{Q}_{4n},\{1,a^2,b\})$ are not isomorphic.
\end{pro}

We first give some definitions and lemmas needed in its proof. The \emph{spectrum} of a graph $\Gamma$ is the spectrum of its adjacency matrix $A$, that is, the multiset of eigenvalues of $A$. We define the \emph{characteristic polynomial} of $\Gamma$ to be the characteristic polynomial of $A$. Let $G$ be a finite group and $V$ be a finite-dimensional vector space over the complex field $\mathbb{C}$. A \emph{$\mathbb{C}$-representation} of $G$ is a group homomorphism $\rho\colon G\to\GL(V)$, and the \emph{degree} of $\rho$ is the dimension of $V$. If no nontrivial subspace $U$ of $V$ is fixed by each $\rho(g)$, then we say $\rho$ is \emph{irreducible}. Let $\rho_1\colon G\to\GL(V_1)$ and $\rho_2\colon G\to\GL(V_2)$ be two $\mathbb{C}$-representations of $G$. We say $\rho_1$ and $\rho_2$ are \emph{equivalent} if there is an isomorphism $h$ from $V_1$ to $V_2$ such that $\rho_1(g)h=h\rho_2(g)$ for each $g\in G$. Otherwise, we say $\rho_1$ and $\rho_2$ are \emph{inequivalent}. Denote by $irr(G)=\{\rho_1,\ldots,\rho_m\}$ the set of all inequivalent irreducible $\mathbb{C}$-representations of $G$. For each $\rho_\ell\in irr(G)$ and each subset $S$ of $G$, we define
\[
  \rho_\ell(S)=\sum_{s\in S}\rho_\ell(s),\ \
   A_\ell(S)=\begin{pmatrix}
    \bm{0}_{d_\ell} & \rho_\ell(S^{-1})\\
    \rho_\ell(S) & \bm{0}_{d_\ell}
    \end{pmatrix},\ \text{ and }\
  \chi_{\rho_\ell(S)}(\lambda)=\det\big(\lambda I_{2d_\ell}-A_\ell(S)\big),
\]
where $d_\ell$ is the degree of $\rho_\ell$, and $\bm{0}_{d_\ell}$ is the $d_\ell\times d_\ell$ matrix with all entries $0$.

For a bi-Cayley graph $\mathrm{BCay}(G,S)$ of a finite group $G$, we have the following result.

\begin{lem}\label{lm:char}
Let $\mathrm{BCay}(G,S)$ be a bi-Cayley graph for some $S\subseteq G$, and let $irr(G)=\{\rho_1,\ldots,\rho_m\}$. Then for each $\ell\in\{1,\ldots, m\}$,
\[
  \chi_{\rho_\ell(S)}(\lambda)=\det\big(\lambda^2 I_{d_\ell}-\rho_\ell(S)\rho_\ell(S^{-1})\big).
\]
Furthermore, the characteristic polynomial of $\mathrm{BCay}(G,S)$ is
\[
  \prod_{\ell=1}^m\big(\chi_{\rho_\ell(S)}(\lambda)\big)^{d_\ell}.
\]
\end{lem}

\begin{proof}
A straightforward calculation gives that
\begin{align*}
\chi_{\rho_\ell(S)}(\lambda)
=&\det\big(\lambda I_{2d_\ell}-A_\ell(S)\big)\\
=&\left| \begin{matrix}
\lambda I_{d_\ell} & -\rho_\ell(S^{-1})\\
-\rho_\ell(S) & \lambda I_{d_\ell}
\end{matrix}\right|\\
=&\left| \begin{matrix}
\lambda I_{d_\ell} & -\rho_\ell(S^{-1})\\
\bm{0}_{d_\ell} & \lambda I_{d_\ell}-\rho_\ell(S)\rho_\ell(S^{-1})/\lambda
\end{matrix}\right| \quad\quad\big(\text{row}\;2\to \text{row}\;2+\rho_\ell(S)/\lambda\big)\\
=&\det(\lambda I_{d_\ell})\cdot \det\Big(\lambda I_{d_\ell}-\rho_\ell(S)\rho_\ell(S^{-1})/\lambda\Big)\\
=&\det(\lambda^2 I_{d_\ell}-\rho_\ell(S)\rho_\ell(S^{-1})),
\end{align*}
which proves the first statement. The second statement follows from~\cite[Theorem~1.1(2)]{AT2}.
\end{proof}

The following Lemma lists all inequivalent irreducible $\mathbb{C}$-representations of $\mathrm{Q}_{4n}$.

\begin{lem}\label{lm:rep}\cite[Lemma~2.3]{HuangLi}
All inequivalent irreducible $\mathbb{C}$-representations of $\mathrm{Q}_{4n}$ are characterized in Table~\ref{table:n odd} when $n$ is odd and in Table~\ref{table:n even} otherwise, where $i$ is the imaginary unit and $\omega=e^{\pi i/n}$ is the $(2n)$-th root of unity. (The first column of both tables are notation of representations, and the last two columns list the image of $a^k$ and $ba^k$ under the  representation in the corresponding row.)
\end{lem}

\setlength{\arrayrulewidth}{0.3mm}
\setlength{\tabcolsep}{14pt}
\renewcommand{\arraystretch}{1.2}

\begin{table}[h]
    \centering
    \caption{Inequivalent irreducible $\mathbb{C}$-representation table of $\mathrm{Q}_{4n}$ when $n$ is odd}
    \begin{tabular}[c]{lll}
        \hline
        Representation & $a^k$ ($0\leq k\leq 2n-1$) & $ba^k$ ($0\leq k\leq 2n-1$)\\
        \hline
        $\psi_1$ & $1$ & $1$ \\
        $\psi_2$ & $1$ & $-1$ \\
        $\psi_3$ & $(-1)^k$ & $(-1)^ki$ \\
        $\psi_4$ & $(-1)^k$ & $(-1)^{k+1}i$ \\
        $\phi_j$ ($1\leq j\leq n-1$, $j$ is odd)
        &$\begin{pmatrix}
            \omega^{kj} & 0\\
            0 & \omega^{-kj}
          \end{pmatrix}$
        &$\begin{pmatrix}
            0 & \omega^{-kj}\\
            -\omega^{kj} & 0
          \end{pmatrix}$
        \\
        $\eta_h$ $\big(1\leq h\leq (n-2)/2\big)$
        &$\begin{pmatrix}
            \omega^{2kh} & 0\\
            0 & \omega^{-2kh}
          \end{pmatrix}$
        &$\begin{pmatrix}
            0 & \omega^{-2kh}\\
            \omega^{2kh} & 0
          \end{pmatrix}$\\[15pt]
        \hline \\
    \end{tabular}
    \label{table:n odd}
\end{table}

\begin{table}[h]
    \centering
    \caption{Inequivalent irreducible $\mathbb{C}$-representation table of $\mathrm{Q}_{4n}$ when $n$ is even}
    \begin{tabular}[c]{lll}
        \hline
        Representation & $a^k$ ($0\leq k\leq 2n-1$) & $ba^k$ ($0\leq k\leq 2n-1$)\\
        \hline
        $\psi_1$ & $1$ & $1$ \\
        $\psi_2$ & $1$ & $-1$ \\
        $\psi_3$ & $(-1)^k$ & $(-1)^k$ \\
        $\psi_4$ & $(-1)^k$ & $(-1)^{k+1}$ \\
        $\phi_j$ ($1\leq j\leq n-1$, $j$ is odd)
        &$\begin{pmatrix}
            \omega^{kj} & 0\\
            0 & \omega^{-kj}
          \end{pmatrix}$
        &$\begin{pmatrix}
            0 & \omega^{-kj}\\
            -\omega^{kj} & 0
          \end{pmatrix}$
        \\
        $\eta_h$ $\big(1\leq h\leq (n-2)/2\big)$
        &$\begin{pmatrix}
            \omega^{2kh} & 0\\
            0 & \omega^{-2kh}
          \end{pmatrix}$
        &$\begin{pmatrix}
            0 & \omega^{-2kh}\\
            \omega^{2kh} & 0
          \end{pmatrix}$\\[15pt]
        \hline \\
    \end{tabular}
    \label{table:n even}
\end{table}

\begin{proof}[Proof of Proposition~\ref{prop:not iso}]
Let $S=\{1,a,b\}$, $T=\{1,a^2,b\}$. We prove the conclusion by showing that the spectrums of $\Gamma_S:=\mathrm{BCay}(G,S)$ and $\Gamma_T:=\mathrm{BCay}(G,T)$ are distinct. Using the notation in Tables~\ref{table:n odd} and~\ref{table:n even}, we obtain that whenever $n$ is odd or even,
\[
irr(Q_{4n})=\{\psi_1,\ldots,\psi_4\}\sqcup\{\phi_j\mid 1\leq j\leq n-1, j\;\text{is odd}\}\sqcup \{\eta_h \mid 1\leq h\leq \lfloor(n-1)/2\rfloor\},
\]
where $\sqcup$ represents the disjoint union. Then it follows from Lemma~\ref{lm:char} and Tables~\ref{table:n odd} and~\ref{table:n even} that
\begin{align*}
\chi_{\psi_1(S)}(\lambda)
  =&\det\left(\lambda^2-\psi_1(S)\psi_1(S^{-1})\right)\\
  =&\lambda^2-\psi_1(\{1,a,b\})\psi_1(\{1,a^{-1},ba^n\})\\
  =&\lambda^2-\big(\psi_1(1)+\psi_1(a)+\psi_1(b)\big)\big(\psi_1(1)+\psi_1(a^{-1})+\psi_1(ba^n)\big)\\
  =&\lambda^2-9.
\end{align*}
Similarly,
\begin{align*}
  \chi_{\psi_2(S)}(\lambda)
  =&\det\left(\lambda^2-\psi_2(S)\psi_2(S^{-1})\right)\\
  =&\lambda^2-\big(\psi_2(1)+\psi_2(a)+\psi_2(b)\big)\big(\psi_2(1)+\psi_2(a^{-1})+\psi_2(ba^n)\big)\\
  =&\lambda^2-1,
\end{align*}
\begin{align*}
  \chi_{\psi_3(S)}(\lambda)
  =&\det\left(\lambda^2-\psi_3(S)\psi_3(S^{-1})\right)\\
  =&\lambda^2-\big(\psi_3(1)+\psi_3(a)+\psi_3(b)\big)\big(\psi_3(1)+\psi_3(a^{-1})+\psi_3(ba^n)\big)\\
  =&\begin{cases}
    \lambda^2-(1-1+i)(1-1-i) &\quad\text{if}\;n\;\text{is odd}\\
    \lambda^2-(1-1+1)(1-1+1) &\quad\text{if}\;n\;\text{is even}
  \end{cases}\\
  =&\lambda^2-1,
\end{align*}
\begin{align*}
  \chi_{\psi_4(S)}(\lambda)
  =&\det\left(\lambda^2-\psi_4(S)\psi_4(S^{-1})\right)\\
  =&\lambda^2-\big(\psi_4(1)+\psi_4(a)+\psi_4(b)\big)\big(\psi_4(1)+\psi_4(a^{-1})+\psi_4(ba^n)\big)\\
  =&\begin{cases}
    \lambda^2-(1-1-i)(1-1+i) &\quad\text{if}\;n\;\text{is odd}\\
    \lambda^2-(1-1-1)(1-1-1) &\quad\text{if}\;n\;\text{is even}
  \end{cases}\\
  =&\lambda^2-1.
\end{align*}
Note that $\omega=e^{\pi i/n}$, and thus $\omega^n=-1$.
For $j\in \{1,\ldots, n-1\}$ and $j$ is odd, we have
\begin{align*}
  \phi_j(S)\phi_j(S^{-1})
  =&\big(\phi_j(1)+\phi_j(a)+\phi_j(b)\big)\big(\phi_j(1)+\phi_j(a^{-1})+\phi_j(ba^n)\big)\\
  =&\left(\begin{pmatrix}
    1 & 0\\
    0 & 1
  \end{pmatrix}
  +\begin{pmatrix}
    \omega^j & 0\\
    0 & \omega^{-j}
  \end{pmatrix}
  +\begin{pmatrix}
    0 & 1 \\
    -1 & 0
  \end{pmatrix}\right)\cdot\\
  &\left(\begin{pmatrix}
    1 & 0\\
    0 & 1
  \end{pmatrix}
  +\begin{pmatrix}
    \omega^{-j} & 0\\
    0 & \omega^j
  \end{pmatrix}
  +\begin{pmatrix}
    0 & \omega^{-nj} \\
    -\omega^{nj} & 0
  \end{pmatrix}\right)\\
  =&\begin{pmatrix}
    1+\omega^j & 1\\
    -1 & 1+\omega^{-j}
  \end{pmatrix}
  \begin{pmatrix}
    1+\omega^{-j} & -1\\
    1 & 1+\omega^{j}
  \end{pmatrix}\\
  =&\begin{pmatrix}
    (1+\omega^j)(1+\omega^{-j})+1 & 0\\
    0 & (1+\omega^j)(1+\omega^{-j})+1
  \end{pmatrix}\\
  =& \big(3+2\cos\frac{j\pi}{n}\big)I_2.
\end{align*}
It follows that
\begin{align*}
  \chi_{\phi_j(S)}(\lambda)
  =&\det\left(\lambda^2I_2-\phi_j(S)\phi_j(S^{-1})\right)\\
  =&\det\left(\lambda^2I_2-\big(3+2\cos\frac{j\pi}{n}\big)I_2\right)\\
  =&\big(\lambda^2-\big(3+2\cos\frac{j\pi}{n}\big)\big)^2.
\end{align*}
Similarly, for $h\in \{1,\ldots,\lfloor (n-1)/2\rfloor\}$, we have
\begin{align*}
  \eta_h(S)\eta_h(S^{-1})
  =&\big(\eta_h(1)+\eta_h(a)+\eta_h(b)\big)\big(\eta_h(1)+\eta_h(a^{-1})+\eta_h(ba^n)\big)\\
  =&\left(\begin{pmatrix}
    1 & 0\\
    0 & 1
  \end{pmatrix}
  +\begin{pmatrix}
    \omega^{2h} & 0\\
    0 & \omega^{-2h}
  \end{pmatrix}
  +\begin{pmatrix}
    0 & 1 \\
    1 & 0
  \end{pmatrix}\right)\cdot\\
  &\left(\begin{pmatrix}
    1 & 0\\
    0 & 1
  \end{pmatrix}
  +\begin{pmatrix}
    \omega^{-2h} & 0\\
    0 & \omega^{2h}
  \end{pmatrix}
  +\begin{pmatrix}
    0 & \omega^{-2nh} \\
    \omega^{2nh} & 0
  \end{pmatrix}\right)\\
  =&\begin{pmatrix}
    1+\omega^{2h} & 1\\
    1 & 1+\omega^{-2h}
  \end{pmatrix}
  \begin{pmatrix}
    1+\omega^{-2h} & 1\\
    1 & 1+\omega^{2h}
  \end{pmatrix}\\
  =&\begin{pmatrix}
    (1+\omega^{2h})(1+\omega^{-2h})+1 & 2(1+\omega^{2h})\\
    2(1+\omega^{-2h}) & (1+\omega^{2h})(1+\omega^{-2h})+1
  \end{pmatrix}.
\end{align*}
We then deduce that
\begin{align*}
  \chi_{\eta_h(S)}(\lambda)
  =&\det\left(\lambda^2I_2-\eta_h(S)\eta_h(S^{-1})\right)\\
  =&\det\left(
    \begin{pmatrix}
      \lambda^2 & 0\\
      0 & \lambda^2
    \end{pmatrix}
    -\begin{pmatrix}
      (1+\omega^{2h})(1+\omega^{-2h})+1 & 2(1+\omega^{2h})\\
      2(1+\omega^{-2h}) & (1+\omega^{2h})(1+\omega^{-2h})+1
    \end{pmatrix}
  \right)\\
  =&\big(\lambda^2-(1+\omega^{2h})(1+\omega^{-2h})-1\big)^2-4(1+\omega^{2h})(1+\omega^{-2h})\\
  =&\big(\lambda^2-\big(3+2\cos\frac{2h\pi}{n}\big)\big)^2-4\big(2+2\cos\frac{2h\pi}{n}\big)\\
  =&\big(\lambda^2-\big(3+2\cos\frac{2h\pi}{n}\big)\big)^2-\big(4\cos\frac{h\pi}{n}\big)^2.
\end{align*}
We obtain from Lemma~\ref{lm:char} that the characteristic polynomial of $\Gamma_S$ is
\[
  \prod_{k=1}^4\chi_{\psi_k(S)}(\lambda)\prod_{\substack{j=1\\ j\;\text{odd}}}^{n-1}\big(\chi_{\phi_j(S)}(\lambda)\big)^2\prod_{h=1}^{\lfloor (n-1)/2\rfloor}\big(\chi_{\eta_h(S)}(\lambda)\big)^2.
\]

Next, we consider the characteristic polynomial of $\Gamma_T$. Along the same way as the above calculations, we obtain
\begin{align*}
\chi_{\psi_1(T)}(\lambda)
  =&\det\left(\lambda^2-\psi_1(T)\psi_1(T^{-1})\right)\\
  =&\lambda^2-\psi_1(\{1,a^2,b\})\psi_1(\{1,a^{-2},ba^n\})\\
  =&\lambda^2-\big(\psi_1(1)+\psi_1(a^2)+\psi_1(b)\big)\big(\psi_1(1)+\psi_1(a^{-2})+\psi_1(ba^n)\big)\\
  =&\lambda^2-9,
\end{align*}
\begin{align*}
  \chi_{\psi_2(T)}(\lambda)
  =&\det\left(\lambda^2-\psi_2(T)\psi_2(T^{-1})\right)\\
  =&\lambda^2-\big(\psi_2(1)+\psi_2(a^2)+\psi_2(b)\big)\big(\psi_2(1)+\psi_2(a^{-2})+\psi_2(ba^n)\big)\\
  =&\lambda^2-1,
\end{align*}
\begin{align*}
  \chi_{\psi_3(T)}(\lambda)
  =&\det\left(\lambda^2-\psi_3(T)\psi_3(T^{-1})\right)\\
  =&\lambda^2-\big(\psi_3(1)+\psi_3(a^2)+\psi_3(b)\big)\big(\psi_3(1)+\psi_3(a^{-2})+\psi_3(ba^n)\big)\\
  =&\begin{cases}
    \lambda^2-(1+1+i)(1+1-i)=\lambda^2-5 &\quad\text{if}\;n\;\text{is odd}\\
    \lambda^2-(1+1+1)(1+1+1)=\lambda^2-9 &\quad\text{if}\;n\;\text{is even}
  \end{cases},
\end{align*}
\begin{align*}
  \chi_{\psi_4(T)}(\lambda)
  =&\det\left(\lambda^2-\psi_4(T)\psi_4(T^{-1})\right)\\
  =&\lambda^2-\big(\psi_4(1)+\psi_4(a^2)+\psi_4(b)\big)\big(\psi_4(1)+\psi_4(a^{-2})+\psi_4(ba^n)\big)\\
  =&\begin{cases}
    \lambda^2-(1+1-i)(1+1+i)=\lambda^2-5 &\quad\text{if}\;n\;\text{is odd}\\
    \lambda^2-(1+1-1)(1+1-1)=\lambda^2-1 &\quad\text{if}\;n\;\text{is even}
  \end{cases}.
\end{align*}
Recall that $\omega=e^{\pi i/n}$, and so $\omega^n=-1$. For $j\in \{1,\ldots, n-1\}$ and $j$ is odd, we have
\begin{align*}
  \phi_j(T)\phi_j(T^{-1})
  =&\big(\phi_j(1)+\phi_j(a^2)+\phi_j(b)\big)\big(\phi_j(1)+\phi_j(a^{-2})+\phi_j(ba^n)\big)\\
  =&\left(\begin{pmatrix}
    1 & 0\\
    0 & 1
  \end{pmatrix}
  +\begin{pmatrix}
    \omega^{2j} & 0\\
    0 & \omega^{-2j}
  \end{pmatrix}
  +\begin{pmatrix}
    0 & 1 \\
    -1 & 0
  \end{pmatrix}\right)\cdot\\
  &\left(\begin{pmatrix}
    1 & 0\\
    0 & 1
  \end{pmatrix}
  +\begin{pmatrix}
    \omega^{-2j} & 0\\
    0 & \omega^{2j}
  \end{pmatrix}
  +\begin{pmatrix}
    0 & \omega^{-nj} \\
    -\omega^{nj} & 0
  \end{pmatrix}\right)\\
  =&\begin{pmatrix}
    1+\omega^{2j} & 1\\
    -1 & 1+\omega^{-2j}
  \end{pmatrix}
  \begin{pmatrix}
    1+\omega^{-2j} & -1\\
    1 & 1+\omega^{2j}
  \end{pmatrix}\\
  =&\begin{pmatrix}
    (1+\omega^{2j})(1+\omega^{-2j})+1 & 0\\
    0 & (1+\omega^{2j})(1+\omega^{-2j})+1
  \end{pmatrix}\\
  =& \big(3+2\cos\frac{2j\pi}{n}\big)I_2.
\end{align*}
It follows that
\begin{align*}
  \chi_{\phi_j(T)}(\lambda)
  =&\det\left(\lambda^2I_2-\phi_j(T)\phi_j(T^{-1})\right)\\
  =&\det\big(\lambda^2I_2-\big(3+2\cos\frac{2j\pi}{n}\big)I_2\big)\\
  =&\big(\lambda^2-\big(3+2\cos\frac{2j\pi}{n}\big)\big)^2.
\end{align*}
Similarly, for $h\in \{1,\ldots,\lfloor (n-1)/2\rfloor\}$, we have
\begin{align*}
  \eta_h(T)\eta_h(T^{-1})
  =&\big(\eta_h(1)+\eta_h(a^2)+\eta_h(b)\big)\big(\eta_h(1)+\eta_h(a^{-2})+\eta_h(ba^n)\big)\\
  =&\left(\begin{pmatrix}
    1 & 0\\
    0 & 1
  \end{pmatrix}
  +\begin{pmatrix}
    \omega^{4h} & 0\\
    0 & \omega^{-4h}
  \end{pmatrix}
  +\begin{pmatrix}
    0 & 1 \\
    1 & 0
  \end{pmatrix}\right)\cdot\\
  &\left(\begin{pmatrix}
    1 & 0\\
    0 & 1
  \end{pmatrix}
  +\begin{pmatrix}
    \omega^{-4h} & 0\\
    0 & \omega^{4h}
  \end{pmatrix}
  +\begin{pmatrix}
    0 & \omega^{-2nh} \\
    \omega^{2nh} & 0
  \end{pmatrix}\right)\\
  =&\begin{pmatrix}
    1+\omega^{4h} & 1\\
    1 & 1+\omega^{-4h}
  \end{pmatrix}
  \begin{pmatrix}
    1+\omega^{-4h} & 1\\
    1 & 1+\omega^{4h}
  \end{pmatrix}\\
  =&\begin{pmatrix}
    (1+\omega^{4h})(1+\omega^{-4h})+1 & 2(1+\omega^{4h})\\
    2(1+\omega^{-4h}) & (1+\omega^{4h})(1+\omega^{-4h})+1
  \end{pmatrix}.
\end{align*}
We then deduce that
\begin{align*}
  \chi_{\eta_h(T)}(\lambda)
  =&\det\left(\lambda^2I_2-\eta_h(T)\eta_h(T^{-1})\right)\\
  =&\det\left(
    \begin{pmatrix}
      \lambda^2 & 0\\
      0 & \lambda^2
    \end{pmatrix}
    -\begin{pmatrix}
      (1+\omega^{4h})(1+\omega^{-4h})+1 & 2(1+\omega^{4h})\\
      2(1+\omega^{-4h}) & (1+\omega^{4h})(1+\omega^{-4h})+1
    \end{pmatrix}
  \right)\\
  =&\big(\lambda^2-(1+\omega^{4h})(1+\omega^{-4h})-1\big)^2-4(1+\omega^{4h})(1+\omega^{-4h})\\
  =&\big(\lambda^2-\big(3+2\cos\frac{4h\pi}{n}\big)\big)^2-4\big(2+2\cos\frac{4h\pi}{n}\big)\\
  =&\big(\lambda^2-\big(3+2\cos\frac{4h\pi}{n}\big)\big)^2-\big(4\cos\frac{2h\pi}{n}\big)^2.
\end{align*}
We obtain from Lemma~\ref{lm:char} that the characteristic polynomial of $\Gamma_T$ is
\[
  \prod_{k=1}^4\chi_{\psi_k(T)}(\lambda)\prod_{\substack{j=1\\ j\;\text{odd}}}^{n-1}\big(\chi_{\phi_j(T)}(\lambda)\big)^2\prod_{h=1}^{\lfloor (n-1)/2\rfloor}\big(\chi_{\eta_h(T)}(\lambda)\big)^2.
\]

Firstly, assume that $n$ is odd. We compare the multiplicity of the eigenvalue $1$ in the spectrums of $\Gamma_S$ and $\Gamma_T$. For each $j\in \{1,\ldots,n-1\}$ with $j$ odd, we have
\[
  \chi_{\phi_j(S)}(1)=\big(2+2\cos\frac{j\pi}{n}\big)^2=16\cos^4\frac{j\pi}{2n}\neq 0.
\]
Moreover, for each $h\in \{1,\ldots,\lfloor (n-1)/2\rfloor\}$,
\[
  \chi_{\eta_h(S)}(1)=\big(2+2\cos\frac{2h\pi}{n}\big)^2-\big(4\cos\frac{h\pi}{n}\big)^2=16\big(\cos^2\frac{h\pi}{n}\big)\big(\cos^2\frac{h\pi}{n}-1\big)\neq 0.
\]
Thus, the multiplicity of $1$ in the spectrum of $\Gamma_S$ is $3$. Similarly, for each $j\in \{1,\ldots,n-1\}$ with $j$ odd, since $n$ is odd, we obtain
\[
  \chi_{\phi_j(T)}(1)=\big(2+2\cos\frac{2j\pi}{n}\big)^2=16\cos^4\frac{j\pi}{n}\neq 0.
\]
For each $h\in \{1,\ldots,\lfloor (n-1)/2\rfloor\}$, since $n$ is odd, we have $2h\pi/n\notin\{0,\pi/2,\pi\}$, and thus
\[
  \chi_{\eta_h(T)}(1)=\big(2+2\cos\frac{4h\pi}{n}\big)^2-\big(4\cos\frac{2h\pi}{n}\big)^2=16\cos^2\frac{2h\pi}{n}\big(\cos^2\frac{2h\pi}{n}-1\big)\neq 0.
\]
Hence, the multiplicity of $1$ in the spectrum of $\Gamma_T$ is $1$. Therefore,  $1$ is an eigenvalue of both $\Gamma_S$ and $\Gamma_T$ but with different multiplicities.

Next, assume that $n$ is even. Since $\chi_{\psi_1(T)}(3)=\chi_{\psi_3(T)}(3)=0$, we obtain that $3$ is an eigenvalue of $\Gamma_T$ with multiplicity at least $2$. We now show that $3$ is also an eigenvalue of $\Gamma_S$ but has multiplicity $1$. Obviously, $\chi_{\psi_1(S)}(3)= 0$, $\chi_{\psi_2(S)}(3)\neq 0$, $\chi_{\psi_3(S)}(3)\neq 0$ and $\chi_{\psi_4(S)}(3)\neq 0$. Moreover, for $j\in \{1,\ldots, n-1\}$ with $j$ is odd, we have
\[
  \chi_{\phi_j(S)}(3)=\big(6-2\cos\frac{j\pi}{n}\big)^2\neq 0.
\]
For $h\in \{1,\ldots,\lfloor (n-1)/2\rfloor\}$, we have
\[
  \chi_{\eta_h(S)}(3)=\big(6-2\cos\frac{2h\pi}{n}\big)^2-\big(4\cos\frac{h\pi}{n}\big)^2>4^2-4^2=0.
\]
Thus, $\chi_{\eta_h(S)}(3)\neq 0$ for each $h\in \{1,\ldots,\lfloor (n-1)/2\rfloor\}$. Therefore,  $3$ is an eigenvalue of both $\Gamma_S$ and $\Gamma_T$ but with different multiplicities.

In summary, the spectrums of $\Gamma_S$ and $\Gamma_T$ are distinct, and thus $\Gamma_S\not\cong \Gamma_T$.
\end{proof}

\section{Proof of Theorem~\ref{mainth1}}\label{sec:4}

In this section, we prove Theorem~\ref{mainth1}.  We begin by demonstrating the necessary and sufficient condition for $\mathrm{Q}_{4n}$ to be a $2$-BCI-group when $n \geq 3$.
\begin{lem}\label{gq-2-BCI}
Let $\mathrm{Q}_{4n}$ be as in~\eqref{eq:Q4n} for some integer $n\geq 3$. Then $\mathrm{Q}_{4n}$ is a $2$-BCI-group if and only if $n$ is odd.
\end{lem}
\begin{proof}
First, we prove the necessity. Suppose to the contrary that $n$ is even. Then $|a^{n/2}|=|b|=|b^{-1}|=4$. By Lemma~\ref{lm:Q4n}\eqref{enu:cQ4n}, $\langle a\rangle$ is characteristic in $\mathrm{Q}_{4n}$, which implies that there does not exist $\sigma\in \Aut(\mathrm{Q}_{4n})$ such that $(a^{n/2})^{\sigma}=b$ or $(a^{n/2})^{\sigma}=b^{-1}$. Hence, $\mathrm{Q}_{4n}$ is not a FIF-group, and so we derive from Proposition~\ref{finitegroup-1-2-3-BCI}\eqref{enu:2prop} that $\mathrm{Q}_{4n}$ is not a $2$-BCI-group, a contradiction. Therefore, $n$ is odd.

Next, we prove the sufficiency. Assume that $n$ is odd. Let $\BCay(\mathrm{Q}_{4n},S)$ be a bi-Cayley graph for some $S\subseteq \mathrm{Q}_{4n}$ with $|S|=2$, and let $\BCay(\mathrm{Q}_{4n},T)\cong \BCay(\mathrm{Q}_{4n},S)$. Since Proposition~\ref{finitegroup-1-2-3-BCI}\eqref{enu:1prop} shows that $\mathrm{Q}_{4n}$ is an $1$-BCI-group, we are left to prove that $T=gS^\alpha$ for some $g\in\mathrm{Q}_{4n}$ and some $\alpha\in\Aut(\mathrm{Q}_{4n})$. By Remark~\ref{rmk:2.7}, we may assume that $S=\{1,s\}$ and $T=\{1,t\}$, where $s,t\in \mathrm{Q}_{4n}$ and $s\neq t$. Then it follows from Lemma~\ref{isomorphism-Bcay}\eqref{enu:bBiCay} that
\[
\BCay(\langle s\rangle,S)\cong \BCay(\langle t\rangle,T).
\]
This implies that $|\langle s\rangle|=|\langle t\rangle|$, and so $|s|=|t|$. Noting that $n$ is odd, we derive from Lemma~\ref{lm:Q4n}\eqref{enu:bQ4n} that either $s,t\in \langle a\rangle$ or $s,t\in b\langle a\rangle$.

Assume that $s,t\in \langle a\rangle$. Then we obtain from Lemma~\ref{lm:cyclic} that $\langle s\rangle=\langle t\rangle$, and thus,
\[
  \BCay(\langle s\rangle,S)
  \cong \BCay(\langle t\rangle,T)
  =\BCay(\langle s\rangle,T).
\]
This combined with Proposition~\ref{finitegroup-1-2-3-BCI}\eqref{enu:3prop} implies that there exist $g\in \langle a^i\rangle\subseteq \mathrm{Q}_{4n}$ and $\delta\in \Aut(\langle a^i\rangle)$ such that $T=gS^\delta$. Applying Lemma~\ref{lm:Q4n}\eqref{enu:eQ4n}, we can extend $\delta$ to an automorphism of $\mathrm{Q}_{4n}$, say $\alpha$. Therefore, we obtain $T=gS^\alpha$, as required.

Assume that $s,t\in b\langle a\rangle$. Then we write
\[
  s=ba^i \text{ and }t=ba^j, \text{ where }i\neq j\in \mathbb{Z}_{2n}.
\]
It is straightforward to verify that
\[
  S^{\sigma_{1,j}{\sigma_{1,i}}^{-1}}=\{1,ba^i\}^{\sigma_{1,j}{\sigma_{1,i}}^{-1}}=\{1,ba^j\}=T.
\]
Moreover, Lemma~\ref{lm:Q4n}\eqref{enu:fQ4n} shows that $\alpha:=\sigma_{1,j}{\sigma_{1,i}}^{-1}\in \Aut(\mathrm{Q}_{4n})$. Therefore, $T=S^\alpha$, completing the proof.
\end{proof}

In preparation for proving the necessary and sufficient condition for $\mathrm{Q}_{4n}$ to be a $3$-BCI-group when $n \geq 3$, we present the following three technical lemmas.

\begin{lem}\label{lm:reduction}
  Let $\mathrm{Q}_{4n}$ be as in~\eqref{eq:Q4n} for some integer $n\geq 2$. If $S$ is a subset of $\mathrm{Q}_{4n}$ with $1\in S$, $|S|=3$ and $S\not\subseteq \langle a\rangle$, then there exist $g\in \mathrm{Q}_{4n}$ and $\alpha\in \mathrm{Aut}(\mathrm{Q}_{4n})$ such that $gS^\alpha=\{1,a^i,b\}$ for some $i\in \mathbb{Z}_{2n}\setminus \{0\}$.
\end{lem}
\begin{proof}
  Since $S\not\subseteq \langle a\rangle$, there exists $j\in \mathbb{Z}_{2n}$ such that $ba^j\in S$. Noting that $\sigma_{1,2n-j}\in \mathrm{Aut}(\mathrm{Q}_{4n})$ and $(ba^j)^{\sigma_{1,2n-j}}=b$, we have $\{1,b\}\subseteq S^{\sigma_{1,2n-j}}$. Let $S^{\sigma_{1,2n-j}}=\{1,b,c\}$. If $c\in \langle a\rangle$, then the lemma holds. Now we assume $c=ba^i$ for some $i\in\mathbb{Z}_{2n}\setminus\{0\}$. Then a straightforward calculation gives that
  \[
    (b^{-1}\{1,b,c\})^{\sigma_{1,n}}=\{b^{-1},1,b^{-1}c\}^{\sigma_{1,n}}=\{ba^n,1,a^i\}^{\sigma_{1,n}}=\{b,1,a^i\}.
  \]
  Since $\sigma_{1,n}\in \mathrm{Aut}(\mathrm{Q}_{4n})$, we conclude that $(b^{-1})^{\sigma_{1,n}}S^{\sigma_{1,2n-j}\sigma_{1,n}}=\{1,a^i,b\}$, as desired.
\end{proof}

\begin{lem}\label{lm:3.5}
  Let $i\in \mathbb{Z}_{2n}\setminus\{0\}$ and $d=\gcd(i,2n)$. Then there exists $\alpha\in \mathrm{Aut(\mathrm{Q}_{4n})}$ such that $\{1,a^i,b\}^\alpha=\{1,a^d,b\}$.
\end{lem}

\begin{proof}
  Let $i=dk$. Then $\gcd(k,2n)=1$, which, by B\'ezout's lemma, implies that there exists $r,s\in \mathbb{Z}$ such that $kr+2ns=1$. Clearly, $\gcd(r,2n)=1$ and $r$ can be regarded as an element of $\mathbb{Z}_{2n}^*$. Hence, $\sigma_{r,0}\in \mathrm{Aut}(\mathrm{Q}_{4n})$ and
  \[
    \{1,a^i,b\}^{\sigma_{r,0}}=\{1,a^{ri},b\}=\{1,a^{rdk},b\}=\{1,a^{d(1-2ns)},b\}=\{1,a^d,b\},
  \]
  completing the proof.
\end{proof}

\begin{lem}\label{gq-3-noisomorphic}
Let $\mathrm{Q}_{4n}$ be as in~\eqref{eq:Q4n} for some odd integer $n\geq 3$. Take $S=\{1,a^i,b\}$ and $T=\{1,a^j,b\}$ such that $2|a^i|=|a^j|$, where $i,j\in\mathbb{Z}_{2n}\setminus\{0\}$. Then $\BCay(\mathrm{Q}_{4n},S)\ncong \BCay(\mathrm{Q}_{4n},T)$.
\end{lem}
\begin{proof}
  We derive from Lemma~\ref{lm:3.5} that there exist $\alpha,\beta\in \mathrm{Aut(\mathrm{Q}_{4n})}$ such that
  \[
    S^\alpha=\{1,a^i,b\}^\alpha=\{1,a^{\gcd(i,2n)},b\}\ \text{ and }\ T^\beta=\{1,a^j,b\}^\beta=\{1,a^{\gcd(j,2n)},b\}.
  \]
  Then it follows from Lemma~\ref{isomorphism-Bcay}\eqref{enu:aBiCay} that
  \[
    \BCay(\mathrm{Q}_{4n}, S)\cong \BCay(\mathrm{Q}_{4n}, S^\alpha)\ \text{ and }\ \BCay(\mathrm{Q}_{4n}, T)\cong \BCay(\mathrm{Q}_{4n}, T^\beta).
  \]
  Hence, by Lemma~\ref{isomorphism-Bcay}\eqref{enu:bBiCay},
  \[
    \BCay(\mathrm{Q}_{4n},S)\ncong \BCay(\mathrm{Q}_{4n},T)\ \text{ if and only if }\
    \BCay(\langle S^\alpha\rangle,S^\alpha)\ncong \BCay(\langle T^\beta\rangle,T^\beta).
  \]
  Since $2|a^i|=|a^j|$, we obtain from Lemma~\ref{lm:Q4n}\eqref{enu:bQ4n} that
  \[
    \gcd(i,2n)=2\gcd(j,2n).
  \]
  This implies that $\gcd(j,2n)$ is odd as $n$ is odd, and moreover, $i,j\in \mathbb{Z}_{2n}\setminus \{0,n\}$. As a consequence, we deduce from Lemma~\ref{lm:Q4n}\eqref{enu:dQ4n} that
  \[
    \langle S^\alpha\rangle=\langle a^{\gcd(i,2n)},b \rangle\cong \mathrm{Q}_{4m} \cong \langle a^{\gcd(j,2n)},b \rangle\cong \langle T^\beta\rangle,
  \]
  where $m=|a^{\gcd(i,2n)}|$. Let $c=a^{\gcd(j,2n)}$. Then $c^2=a^{\gcd(i,2n)}$, and so,
  \[
    \BCay(\langle S^\alpha\rangle,S^\alpha)\cong
    \BCay(\mathrm{Q}_{4m}, \{1,c^2,b\})\ \text{ and }\ \BCay(\langle T^\beta\rangle,S^\beta)\cong
    \BCay(\mathrm{Q}_{4m}, \{1,c,b\}).
  \]
  Noting by Proposition~\ref{prop:not iso} that
  \[
    \BCay(\mathrm{Q}_{4m}, \{1,c^2,b\})\not \cong \BCay(\mathrm{Q}_{4m}, \{1,c,b\}),
  \]
  we conclude that $\BCay(\mathrm{Q}_{4n},S)\ncong \BCay(\mathrm{Q}_{4n},T)$, as required.
\end{proof}

Now, we are ready to prove the necessary and sufficient condition for $\mathrm{Q}_{4n}$ to be a $3$-BCI-group when $n \geq 3$, which constitutes the main part of the proof of Theorem~\ref{mainth1}.

\begin{lem}\label{gq-3-BCI}
  Let $\mathrm{Q}_{4n}$ be as in~\eqref{eq:Q4n} for some integer $n\geq 3$. Then $\mathrm{Q}_{4n}$ is a $3$-BCI-group if and only if $n$ is odd.
\end{lem}
\begin{proof}
 The necessity follows from Lemma~\ref{gq-2-BCI}. Now we prove the sufficiency. Assume that $n$ is odd. By Lemma~\ref{gq-2-BCI}, it suffices to prove that for each subset $S\subseteq \mathrm{Q}_{4n}$ with $|S|=3$, the graph $\BCay(\mathrm{Q}_{4n},S)$ is a BCI-graph. Let $\BCay(\mathrm{Q}_{4n},T)\cong \BCay(\mathrm{Q}_{4n},S)$ for some $T\subseteq \mathrm{Q}_{4n}$. Clearly, $|T|=|S|=3$. By Remark~\ref{rmk:2.7}, we may assume that $1\in S$ and $1\in T$. This together with Lemma~\ref{isomorphism-Bcay}\eqref{enu:bBiCay} gives that $\BCay(\langle S\rangle, S)\cong \BCay(\langle T\rangle,T)$, and so, $|\langle S\rangle|=|\langle T\rangle|$. To finish the proof, we prove in the following two cases that $T=gS^\alpha$ for some $g\in\mathrm{Q}_{4n}$ and some $\alpha\in\Aut(\mathrm{Q}_{4n})$.

\smallskip
\noindent{\bf Case 1:} $S\subseteq \langle a\rangle$.
\smallskip

In this case, $|\langle S\rangle|$ divides $2n$. Since $n$ is odd, it follows from $|\langle S\rangle|=|\langle T\rangle|$ that $T\subseteq \langle a\rangle$. Then we derive from Lemma~\ref{lm:cyclic} that $\langle S\rangle=\langle T\rangle$ is a cyclic subgroup of $ \langle a\rangle$. Thus, $\BCay(\langle S\rangle, S)\cong \BCay(\langle T\rangle,T)\cong \BCay(\langle S\rangle,T)$. It follows from Proposition~\ref{finitegroup-1-2-3-BCI}(iii) that there exist $g\in \langle S\rangle\subseteq \mathrm{Q}_{4n}$ and $\delta\in \Aut(\langle S\rangle)$ such that $T=gS^\delta$. By Lemma~\ref{lm:Q4n}\eqref{enu:eQ4n}, $\delta$ can be extended to an automorphism of $\mathrm{Q}_{4n}$, say $\alpha$. Consequently, $T=gS^\alpha$, as required.
\smallskip

\noindent{\bf Case 2:} $S\not\subseteq \langle a\rangle$.
\smallskip

In this case, we have $ba^r\in S$ for some $r\in \mathbb{Z}_{2n}$. Then we see from Lemma~\ref{lm:Q4n}\eqref{enu:bQ4n} that $4$ divides $|\langle S\rangle|=|\langle T\rangle|$, and thus, $T\not\subseteq \langle a\rangle$ as $n$ is odd. Therefore, it follows from Lemma~\ref{lm:reduction} that there exist $g,h\in \mathrm{Q}_{4n}$ and $\alpha,\beta\in \mathrm{Aut}(\mathrm{Q}_{4n})$ such that
\begin{equation}\label{eq:S'T'}
  S':=gS^\alpha=\{1,a^i,b\}\ \text{ and }\ T':=hT^\beta=\{1,a^j,b\}
\end{equation}
for some $i,j\in \mathbb{Z}_{2n}\setminus\{0\}$. If $S'=T'$, then $T=(h^{-1}g)^{\beta^{-1}}S^{\alpha\beta^{-1}}$, as required. Next, we assume that $i\neq j$.
We see from Lemma~\ref{isomorphism-Bcay}\eqref{enu:aBiCay} that
\[
  \BCay(\mathrm{Q}_{4n}, S')\cong \BCay(\mathrm{Q}_{4n}, S) \cong \BCay(\mathrm{Q}_{4n}, T)\cong \BCay(\mathrm{Q}_{4n}, T').
\]
This in conjunction with Lemma~\ref{isomorphism-Bcay}\eqref{enu:bBiCay} gives that
\[
  \BCay(\langle S'\rangle, S')\cong \BCay(\langle T'\rangle, T').
\]
Thus, $|\langle a^i,b\rangle|=|\langle S'\rangle|=|\langle T'\rangle|=|\langle a^j,b\rangle|$. We deduce from Lemma~\ref{lm:Q4n}\eqref{enu:dQ4n} that \[
  |a^i|=|a^j|,\ \text{ or }\ 2|a^i|=|a^j|,\ \text{ or }\ |a^i|=2|a^j|.
\]
For the latter two cases, Lemma~\ref{gq-3-noisomorphic} asserts that $\BCay(\mathrm{Q}_{4n},S)\ncong \BCay(\mathrm{Q}_{4n},T)$, a contradiction. Hence $|a^i|=|a^j|$, and therefore we obtain from Lemma~\ref{lm:Q4n}\eqref{enu:bQ4n} that $\gcd(i,2n)=\gcd(j,2n)$. Then we derive from Lemma~\ref{lm:3.5} that there exist $\gamma,\omega\in \mathrm{Aut(\mathrm{Q}_{4n})}$ such that
\[
  (S')^\gamma=\{1,a^i,b\}^\gamma=\{1,a^{\gcd(i,2n)},b\}=\{1,a^{\gcd(j,2n)},b\}=\{1,a^j,b\}^\omega=(T')^\omega.
\]
Therefore, this combined with~\eqref{eq:S'T'} concludes that $(gS^\alpha)^\gamma=(hT^\beta)^\omega$, and so,
\[
  T=(h^{-1}(gS^\alpha)^{\gamma\omega^{-1}})^{\beta^{-1}}=(h^{-1}g^{\gamma\omega^{-1}})^{\beta^{-1}}S^{\alpha\gamma\omega^{-1}\beta^{-1}}.
\]
This completes the proof.
\end{proof}

\medskip

\begin{proof}[Proof of Theorem~\ref{mainth1}]
  Lemmas~\ref{gq-2-BCI} and~\ref{gq-3-BCI} shows that if $n\geq 3$, then $\mathrm{Q}_{4n}$ is a $3$-BCI-group if and only if $\mathrm{Q}_{4n}$ is a $2$-BCI-group if and only if $n$ is odd. Furthermore, we see from~\cite[Lemma 3.4]{JinL2} that $\mathrm{Q}_8$ is a BCI-group, and thus, $\mathrm{Q}_8$ is an $m$-BCI-group for each $m\in \{2,3\}$. This completes the proof.
\end{proof}


\noindent
{\bf Acknowledgements.} The second author was supported by the Melbourne Research Scholarship provided by The University of Melbourne.

\end{document}